%% file: stocon.tex
\documentclass[a4paper,11pt]{article}

\usepackage{accents}
\usepackage{amssymb}
\usepackage{graphicx}
\usepackage{setspace}
\usepackage{natbib}
\usepackage{mathrsfs}

\newcommand{\bfa}{\mathbf{a}}
\newcommand{\bfb}{\mathbf{b}}
\newcommand{\bfc}{\mathbf{c}}

\newcommand{\bff}{\mathbf{f}}

\newcommand{\bfx}{\mathbf{x}}
\newcommand{\bfy}{\mathbf{y}}

\newcommand{\bfA}{\mathbf{A}}
\newcommand{\bfB}{\mathbf{B}}

\newcommand{\bfF}{\mathbf{F}}

\newcommand{\bfH}{\mathbf{H}}
\newcommand{\bfI}{\mathbf{I}}
\newcommand{\bfJ}{\mathbf{J}}
\newcommand{\bfK}{\mathbf{K}}
\newcommand{\bfL}{\mathbf{L}}
\newcommand{\bfM}{\mathbf{M}}

\newcommand{\bfV}{\mathbf{V}}

\newcommand{\bfTh}{\mathbf{\Theta}}
\newcommand{\zeros}{\mathbf{0}}
\newcommand{\sM}{\mathcal{M}}

\newcommand{\expect}{\mathbb{E}}
\newcommand{\expectx}{\mathbb{E}_{\bfx_0}}
\newcommand{\proba}{\mathbb{P}}
\newcommand{\probax}{\mathbb{P}_{\bfx_0}}

\makeatletter

\@addtoreset{equation}{section}
\makeatother

\newcommand{\bffg}{\accentset{\frown}{\mathbf{f}}}
\newcommand{\bfxg}{\accentset{\frown}{\mathbf{x}}}
\newcommand{\bfyg}{\accentset{\frown}{\mathbf{y}}}
\newcommand{\sg}{\accentset{\frown}{\sigma}}

\newtheorem{lemma}{Lemma}
\newtheorem{theorem}{Theorem}
\newtheorem{cor}{Corollary}
\newtheorem{defi}{Definition}

\title{A Contraction Theory Approach to Stochastic Incremental
Stability}

\author{Quang-Cuong Pham\,\thanks{To whom correspondance should be addressed.}\\
{\normalsize LPPA, Coll\`ege de France}\\
{\normalsize Paris, France}\\
{\normalsize\texttt{cuong.pham@normalesup.org}}
\and Nicolas Tabareau\\
{\normalsize LPPA, Coll\`ege de France}\\
{\normalsize Paris, France}\\
{\normalsize\texttt{tabareau.nicolas@gmail.com}}
\and Jean-Jacques Slotine\\
{\normalsize Nonlinear Systems Laboratory, MIT}\\
{\normalsize Cambridge, MA 02139, USA}\\
{\normalsize\texttt{jjs@mit.edu}}}

\begin{document}

\maketitle

\begin{abstract}
  We investigate the incremental stability properties of It\^o
  stochastic dynamical systems. Specifically, we derive a stochastic
  version of nonlinear contraction theory that provides a bound on the
  mean square distance between any two trajectories of a
  stochastically contracting system. This bound can be expressed as a
  function of the noise intensity and the contraction rate of the
  noise-free system. We illustrate these results in the contexts of
  stochastic nonlinear observers design and stochastic
  synchronization.
\end{abstract}

\input{intro.tex}

\input{result.tex}

\input{combinations.tex}

\input{examples.tex}

\paragraph{Acknowledgments} We are grateful to Dr S. Darses, Prof D.
Bennequin and Prof M. Yor for stimulating discussions, and to the
Associate Editor and the reviewers for their helpful comments.

\input{appendix.tex}

\bibliographystyle{abbrv}
\bibliography{nsl}

\end{document}

%% file: intro.tex
\section{Introduction}

Nonlinear stability properties are often considered with respect to an
equilibrium point or to a nominal system trajectory (see e.g.
\cite{SL91}). By contrast, \emph{incremental} stability is concerned
with the behaviour of system trajectories \emph{with respect to each
  other}. From the triangle inequality, global exponential incremental
stability (any two trajectories tend to each other exponentially) is a
stronger property than global exponential convergence to a single
trajectory.

Historically, work on deterministic incremental stability can be
traced back to the 1950's \cite{Lew49,Dem61,Har64} (see e.g.
\cite{LS05,Jou05} for a more extensive list and historical discussion
of related references). More recently, and largely independently of
these earlier studies, a number of works have put incremental
stability on a broader theoretical basis and made relations with more
traditional stability approaches \cite{Fro97,SW97,LS98,Ang02,AC05}.
Furthermore, it was shown that incremental stability is especially
relevant in the study of such problems as state detection
\cite{Ang02}, observer design or synchronization analysis.

While the above references are mostly concerned with
\emph{deterministic} stability notions, stability theory has also been
extended to \emph{stochastic} dynamical systems, see for
instance~\cite{Kus67,Has80}. This includes important recent
developments in Lyapunov-like approaches~\cite{Flo95,Mao91}, as
well as applications to standard problems in systems and control
\cite{Flo97,Tsi99,DenX01}. However, stochastic versions
of incremental stability have not yet been systematically
investigated.

The goal of this paper is to extend some concepts and results in
incremental stability to stochastic dynamical systems. More
specifically, we derive a stochastic version of contraction analysis
in the specialized context of state-independent metrics. 

We prove in section~\ref{sec:results} that the mean square distance
between any two trajectories of a stochastically contracting system is
upper-bounded by a constant after exponential transients. In contrast
with previous works on incremental stochastic stability \cite{CarX04},
we consider the case when the two trajectories are subject to
\emph{distinct} and independent noises, as detailed in
section~\ref{sec:settings}. This specificity enables our theory to
have a number of new and practically important applications. However,
the fact that the noise does not vanish as two trajectories get very
close to each other will prevent us from obtaining asymptotic
almost-sure stability results (see section~\ref{sec:noalmostsure}).

In section~\ref{sec:combinations}, we show that results on
combinations of deterministic contracting systems have simple
analogues in the stochastic case. These combination properties allow
one to build by recursion stochastically contracting systems of
arbitrary size.

Finally, as illustrations of our results, we study in section
\ref{sec:examples} several examples, including contracting observers
with noisy measurements, stochastic composite variables and
synchronization phenomena in networks of noisy dynamical systems.


%% file: result.tex
\section{Main results}
\label{sec:results}

\subsection{Background}

\subsubsection{Nonlinear contraction theory}

Contraction theory~\cite{LS98} provides a set of tools to analyze the
incremental exponential stability of nonlinear systems, and has been
applied notably to observer design~\cite{LS98,LS00,AR03,JF04,ZS05},
synchronization analysis~\cite{WS05,PS07} and systems neuroscience
modelling~\cite{GirX08}. Nonlinear contracting systems enjoy desirable
aggregation properties, in that contraction is preserved under many
types of system combinations given suitable simple
conditions~\cite{LS98}.

While we shall derive global properties of nonlinear systems, many of
our results can be expressed in terms of eigenvalues of symmetric
matrices~\cite{HJ85}. Given a square matrix $\bfA$, the
symmetric part of $\bfA$ is denoted by $\bfA_s$. The smallest and
largest eigenvalues of $\bfA_s$ are denoted by
$\lambda_\mathrm{min}(\bfA)$ and $\lambda_\mathrm{max}(\bfA)$. Given
these notations, the matrix $\bfA$ is {\it positive definite} (denoted
$\bfA > \zeros$) if $\lambda_\mathrm{min}(\bfA)>0$, and it is
\emph{uniformly} positive definite if 
\[
\exists \beta >0 \quad  \forall
  \bfx,t \quad \lambda_\mathrm{min}(\bfA(\bfx,t)) \ge \beta
\]

The basic theorem of contraction analysis, derived in~\cite{LS98},
can be stated as follows

\begin{theorem}[Contraction]
\label{theorem:contraction}
Consider, in $\mathbb{R}^n$, the deterministic system
\begin{equation}
  \label{equ:main2}
  \dot\bfx = \bff(\bfx,t)
\end{equation}
where $\bff$ is a smooth nonlinear function. Denote the Jacobian
matrix of $\bff$ with respect to its first variable by $\frac{\partial
  \bff} {\partial\bfx}$. If there exists a square matrix
$\bfTh(\bfx,t)$ such that
$\bfM(\bfx,t)=\bfTh(\bfx,t)^T\bfTh(\bfx,t)$ is uniformly positive
definite and the matrix
\[
\bfF(\bfx,t) = \left(\frac{d}{dt}\bfTh(\bfx,t) + \bfTh(\bfx,t)
  \frac{\partial \bff} {\partial\bfx} \right) \bfTh^{-1}(\bfx,t)
\]
is uniformly negative definite, then all system trajectories converge
exponentially to a single trajectory, with convergence rate
$|\sup_{\bfx,t}\lambda_\mathrm{max}(\bfF)|=\lambda >0$. The system is
said to be \emph{contracting}, $\bfF$ is called its \emph{generalized
  Jacobian}, $\bfM(\bfx,t)$ its contraction \emph{metric} and
$\lambda$ its contraction \emph{rate}.
\end{theorem}

\subsubsection{Standard stochastic stability}

\label{sec:standard-stab}

In this section, we present very informally the basic ideas of
standard stochastic stability (for a rigourous treatment, the reader
is referred to e.g. \cite{Kus67}). This will set the context to
understand the forthcoming difficulties and differences associated
with incremental stochastic stability.

For simplicity, we consider the special case of global exponential
stability. Let $\bfx(t)$ be a Markov stochastic process and assume
that there exists a non-negative function $V$ ($V(\bfx)$ may represent
e.g. the squared distance of $\bfx$ from the origin) such that
\begin{equation}
  \label{eq:ideal-bound}
  \forall \bfx\in\mathbb{R}^n \quad \widetilde{A}V(\bfx) \leq -\lambda
  V(\bfx)
\end{equation}
where $\lambda$ is a positive real number and $\widetilde{A}$ is the
infinitesimal operator of the process $\bfx(t)$. The operator
$\widetilde{A}$ is the stochastic analogue of the deterministic
differentiation operator. In the case that $\bfx(t)$ is an It\^o
process, $\widetilde{A}$ corresponds to the widely-used
\cite{Mao91,Tsi99,DenX01} differential generator $\mathscr{L}$ (for a
proof of this fact, see \cite{Kus67}, p.~15 or \cite{Arn74}, p.~42).

For $\bfx_0\in\mathbb{R}^n$, let
$\expectx(\cdot)=\expect(\cdot|\bfx(0)=\bfx_0)$. Then by
Dynkin's formula (\cite{Kus67}, p.~10), one has
\[
\begin{array}{rcl}
  \forall t\geq 0 \quad \expectx V(\bfx(t))-V(\bfx_0) &=& \expectx \int_0^t
  \widetilde{A}V(\bfx(s)) ds \\
  &\leq& -\lambda \expectx \int_0^t V(\bfx(s))ds 
  = -\lambda  \int_0^t \expectx V(\bfx(s))ds 
\end{array}
\]

Applying the Gronwall's lemma to the deterministic real-valued
function $t\to\expectx V(\bfx(t))$ yields
\[
\forall t\geq 0 \quad \expectx V(\bfx(t)) \leq  V(\bfx_0) e^{-\lambda t}
\]

If we assume furthermore that \mbox{$\expectx V(\bfx(t))<\infty$} for
all $t$, then the above implies that $V(\bfx(t))$ is a supermartingale
(see lemma~\ref{lemma:super} in the Appendix for details), which
yields, by the supermartingale inequality
\begin{equation}
  \label{eq:super}
  \probax \left(\sup_{T\leq t < \infty } V(\bfx(t)) \geq A \right) \leq \frac{
    \expectx V(\bfx(T))}{A} \leq \frac{V(\bfx_0) e^{-\lambda T}}{A}
\end{equation}

Thus, one obtains an \emph{almost-sure} stability result, in the sense
that
\begin{equation}
  \label{eq:almost-sure}
  \forall A>0 \quad  \lim_{T\to \infty} \probax \left(\sup_{T\leq t <
      \infty } V(\bfx(t)) \geq A \right) =0
\end{equation}

\subsection{The stochastic contraction theorem}

\subsubsection{Settings}
\label{sec:settings}

Consider a noisy system described by an It\^o stochastic differential
equation
\begin{equation}
  \label{eq:main}
  d\bfa=\bff(\bfa,t)dt+\sigma(\bfa,t)dW^d
\end{equation}
where $\bff$ is a $\mathbb{R}^n \times \mathbb{R}^+\to\mathbb{R}^n$
function, $\sigma$ is a $\mathbb{R}^n \times
\mathbb{R}^+\to\mathbb{R}^{nd}$ matrix-valued function and $W^d$ is a
standard $d$-dimensional Wiener process.

To ensure existence and uniqueness of solutions to equation
(\ref{eq:main}), we assume, here and in the remainder of the paper,
the following standard conditions on $\bff$ and $\sigma$

\emph{Lipschitz condition:} There exists a constant $K_1>0$ such that
\[
\forall t\geq 0, \ \bfa,\bfb \in\mathbb{R}^n \quad
\|\bff(\bfa,t)-\bff(\bfb,t)|+\|\sigma(\bfa,t)-\sigma(\bfb,t)\| \leq
K_1\|\bfa-\bfb\|
\]

\emph{Restriction on growth:} There exists a constant $K_2>0$
\[
\forall t\geq 0, \ \bfa \in\mathbb{R}^n \quad
\|\bff(\bfa,t)\|^2+\|\sigma(\bfa,t)\|^2 \leq K_2(1+\|\bfa\|^2)
\]

Under these conditions, one can show (\cite{Arn74}, p.~105) that
equation (\ref{eq:main}) has on $[0,\infty[$ a unique
$\mathbb{R}^n$-valued solution $\bfa(t)$, continuous with probability
one, and satisfying the initial condition $\bfa(0)=\bfa_0$, with
$\bfa_0\in\mathbb{R}^n$. 

In order to investigate the incremental stability properties of
system~(\ref{eq:main}), consider now two system trajectories $\bfa(t)$
and $\bfb(t)$. Our goal will consist of studying the trajectories
$\bfa(t)$ and $\bfb(t)$ with respect to each other. For this, we
consider the \emph{augmented} system $\bfx(t)= (\bfa(t), \bfb(t))^T$,
which follows the equation
\begin{eqnarray}
  \label{eq:dup}
  d\bfx&=&
  \left(\begin{array}{l}
      \bff(\bfa,t)\\
      \bff(\bfb,t)
    \end{array}\right)dt+
  \left(\begin{array}{cc}
      \sigma(\bfa,t)&0\\
      0&\sigma(\bfb,t)
    \end{array}\right)
  \left(\begin{array}{c}
      dW^d_1\\
      dW^d_2
    \end{array}\right)\nonumber\\
  &=&\bffg(\bfx,t)dt+\sg(\bfx,t)dW^{2d}
\end{eqnarray}

\paragraph{Important remark} As stated in the introduction, the
systems $\bfa$ and $\bfb$ are driven by \emph{distinct} and
independent Wiener processes $W_1^d$ and $W_2^d$. This makes our
approach considerably different from \cite{CarX04}, where the authors
studied two trajectories driven by \emph{the same} Wiener process. 

Our approach enables us to study the stability of the system with
respect to variations in initial conditions \emph{and} to random
perturbations: indeed, two trajectories of any real-life system are
typically affected by distinct ``realizations'' of the noise. In
addition, it leads very naturally to nice results on the comparison of
noisy and noise-free trajectories (cf. section~\ref{sec:unperturbed}),
which are particularly useful in applications (cf.
section~\ref{sec:examples}).

However, because of the very fact that the two trajectories are driven
by distinct Wiener processes, we cannot expect the influence of the
noise to vanish when the two trajectories get very close to each
other. This constrasts with~\cite{CarX04}, and more generally, with
the standard stochastic stability case, where the noise vanishes near
the origin (cf. section~\ref{sec:standard-stab}). The consequences of
this will be discussed in detail in section~\ref{sec:noalmostsure}.

\subsubsection{The basic stochastic contraction theorem}
\label{sec:stocon-theorem}

We introduce two hypotheses
\begin{description}
\item[(H1)] $\bff(\bfa,t)$ is contracting in the identity metric, with
  contraction rate $\lambda$, (i.e.  $\forall \bfa,t \quad
  \lambda_{\max}\left(\frac{\partial \bff}{\partial \bfa}\right) \leq
  -\lambda$)
\item[(H2)] $\mathrm{tr}\left(\sigma(\bfa,t)^T\sigma(\bfa,t)\right)$
  is uniformly upper-bounded by a constant $C$ (i.e.  $\forall \bfa,t
  \quad \mathrm{tr}\left(\sigma(\bfa,t)^T\sigma(\bfa,t)\right)\leq C$)
\end{description}

In other words, \textbf{(H1)} says that the noise-free system is
\emph{contracting}, while \textbf{(H2)} says that the variance of the
noise is upper-bounded by a constant.

\begin{defi}
  A system that verifies \emph{\textbf{(H1)}} and \emph{\textbf{(H2)}}
  is said to be \emph{stochastically contracting} in the identity
  metric, with rate $\lambda$ and bound~$C$.
\end{defi}

Consider now the Lyapunov-like function
$V(\bfx)=\|\bfa-\bfb\|^2=(\bfa-\bfb)^T(\bfa-\bfb)$. Using
\textbf{(H1)} and \textbf{(H2)}, we derive below an inequality on
$\widetilde{A}V(\bfx)$, similar to equation (\ref{eq:ideal-bound}) in
section \ref{sec:standard-stab}.

\begin{lemma}
  \label{lemma:ineq}
  Under \emph{\textbf{(H1)}} and \emph{\textbf{(H2)}}, one has the
  inequality
  \begin{equation}
    \label{eq:incremental-bound}
    \widetilde{A}V(\bfx)\leq -2\lambda V(\bfx)+2C
  \end{equation}
\end{lemma}

\paragraph{Proof}

Since $\bfx(t)$ is an It\^o process, $\widetilde{A}$ is given by the
differential operator $\mathscr{L}$ of the process \cite{Kus67,Arn74}.
Thus, by the It\^o formula

\begin{eqnarray}
\widetilde{A}V(\bfx)=\mathscr{L}V(\bfx)&=&\frac{\partial
  V(\bfx)}{\partial \bfx}\bffg(\bfx,t) +\frac{1}{2}\mathrm{tr}\left(
  \sg(\bfx,t)^T\frac{\partial^2 V(\bfx)}{\partial
      \bfx^2}\sg(\bfx,t)\right)
   \nonumber\\ 
&=&\sum_{1\leq i\leq 2n}\frac{\partial V}{\partial \bfx_i} \bffg(\bfx,t)_i
+ \frac{1}{2}\sum_{1\leq i,j,k \leq 2n}\sg(\bfx,t)_{ij}\frac{\partial^2 V}
  {\partial \bfx_i\partial \bfx_k} \sg(\bfx,t)_{kj}
                    \nonumber\\
&=&\sum_{1\leq i\leq n}\frac{\partial V}{\partial \bfa_i} \bff(\bfa,t)_i+
   \sum_{1\leq i\leq n}\frac{\partial V}{\partial \bfb_i} \bff(\bfb,t)_i\nonumber\\
&& + \frac{1}{2}\sum_{1\leq i,j,k \leq n}\sigma(\bfa,t)_{ij}\frac{\partial^2 V}
  {\partial \bfa_i\partial \bfa_k} \sigma(\bfa,t)_{kj}\nonumber\\
&&   + \frac{1}{2}\sum_{1\leq i,j,k \leq n}\sigma(\bfb,t)_{ij}\frac{\partial^2 V}
  {\partial \bfb_i\partial \bfb_k} \sigma(\bfb,t)_{kj} \nonumber\\
&=&2(\bfa-\bfb)^T(\bff(\bfa,t)-\bff(\bfb,t))\nonumber\\
&&   +\mathrm{tr}(\sigma(\bfa,t)^T\sigma(\bfa,t))
   +\mathrm{tr}(\sigma(\bfb,t)^T\sigma(\bfb,t))\nonumber
\end{eqnarray}

Fix $t\geq 0$ and, as in \cite{ES06}, consider the real-valued
function
\[
r(\mu)=(\bfa-\bfb)^T(\bff(\mu \bfa+(1-\mu)\bfb,t)-\bff(\bfb,t))
\]

Since $\bff$ is $C^1$, $r$ is $C^1$ over $[0,1]$. By the mean value
theorem, there exists $\mu_0\in]0,1[$ such that
\[
r'(\mu_0)=r(1)-r(0)=(\bfa-\bfb)^T(\bff(a)-\bff(b))
\]

On the other hand, one obtains by differentiating $r$
\[
r'(\mu_0)=(\bfa-\bfb)^T\left(\frac{\partial \bff}{\partial \bfa}(\mu_0
\bfa+(1-\mu_0)\bfb,t)\right)(\bfa-\bfb)
\]

Thus, one has
\begin{eqnarray}
(\bfa-\bfb)^T(\bff(\bfa)-\bff(\bfb))&=&(\bfa-\bfb)^T\left(\frac{\partial
    \bff}{\partial \bfa}(\mu_0 \bfa+(1-\mu_0)\bfb,t)\right)(\bfa-\bfb)  \nonumber \\
&\leq&-\lambda (\bfa-\bfb)^T(\bfa-\bfb)=-2\lambda V(\bfx)
\end{eqnarray}
where the inequality is obtained by using \textbf{(H1)}.

Finally,
\begin{eqnarray}
\widetilde{A}V(\bfx)&=&2(\bfa-\bfb)^T(\bff(\bfa)-\bff(\bfb))+
\mathrm{tr}(\sigma(\bfa,t)^T\sigma(\bfa,t))+\mathrm{tr}(\sigma(\bfb,t)^T\sigma(\bfb,t))\nonumber \\
&\leq& -2\lambda V(\bfx)+2C \nonumber
\end{eqnarray}
where the inequality is obtained by using \textbf{(H2)}. $\Box$

We are now in a position to prove our main theorem on stochastic
incremental stability.

\begin{theorem}[Stochastic contraction]
  \label{theo:main}
  Assume that system (\ref{eq:main}) verifies \emph{\textbf{(H1)}} and
  \emph{\textbf{(H2)}}. Let $\bfa(t)$ and $\bfb(t)$ be two
  trajectories whose initial conditions are given by a probability
  distribution $p(\bfx(0))=p(\bfa(0),\bfb(0))$. Then
  \begin{equation}
    \label{eq:etot} 
    \forall t \geq 0 \quad \expect \left(\|\bfa(t)-\bfb(t)\|^2 \right) \leq
    \frac{C}{\lambda} +
    e^{-2\lambda t}\int
    \left[\|\bfa_0-\bfb_0\|^2-\frac{C}{\lambda}\right]^+dp(\bfa_0,\bfb_0)
  \end{equation}
  where $[\cdot]^+=\max(0,\cdot)$. This implies in particular
  \begin{equation}
    \label{eq:e} 
    \forall t \geq 0 \quad \expect \left(\|\bfa(t)-\bfb(t)\|^2 \right) \leq
    \frac{C}{\lambda} +
    \expect\left(\|\bfa(0)-\bfb(0)\|^2\right) e^{-2\lambda t}
  \end{equation}
\end{theorem}

\paragraph{Proof}

Let $\bfx_0=(\bfa_0,\bfb_0)\in\mathbb{R}^{2n}$. By Dynkin's formula
(\cite{Kus67}, p.~10)
\[
\expectx V(\bfx(t))-V(\bfx_0)=\expectx\int_0^t \widetilde{A}V(\bfx(s))ds
\]

Thus one has $\forall u,t \quad 0\leq u \leq t < \infty$
\begin{eqnarray}
  \expectx V(\bfx(t))-\expectx V(\bfx(u))&=&\expectx 
  \int_u^t \widetilde{A}V(\bfx(s))ds
  \nonumber\\ 
  \label{eq:a1}
  &\leq&\expectx\int_u^t (-2\lambda V(\bfx(s)) + 2C) ds \\
  &=&\int_u^t (-2\lambda \expectx V(\bfx(s)) + 2C)ds \nonumber
\end{eqnarray}
where inequality (\ref{eq:a1}) is obtained by using lemma
\ref{lemma:ineq}.

Denote by $g(t)$ the \emph{deterministic} quantity $\expectx
V(\bfx(t))$.  Clearly, $g(t)$ is a continuous function of $t$ since
$\bfx(t)$ is a continuous process. The function $g$ then satisfies
the conditions of the Gronwall-type lemma \ref{lemma:gronwall} in the
Appendix, and as a consequence
\[
\forall t \geq 0 \quad \expectx V(\bfx(t)) \leq \frac{C}{\lambda}+
\left[V(\bfx_0)-\frac{C}{\lambda}\right]^+e^{-2\lambda T}
\]

Integrating the above inequality with respect to $\bfx_0$ yields the
desired result (\ref{eq:etot}). Next, inequality (\ref{eq:e}) follows
from (\ref{eq:etot}) by remarking that
\begin{eqnarray}
  \label{eq:trick}
  \int
  \left[\|\bfa_0-\bfb_0\|^2-\frac{C}{\lambda}\right]^+dp(\bfa_0,\bfb_0) &\leq&
  \int 
  \|\bfa_0-\bfb_0\|^2 dp(\bfa_0,\bfb_0) \nonumber \\
  &=& \expect\left(\|\bfa(0)-\bfb(0)\|^2\right) \quad 
\end{eqnarray} $\Box$

\paragraph{Remark} Let $\epsilon>0$ and $T_\epsilon=\frac{1}{2\lambda}
\log\left(\sqrt{\frac{\expect\left(\|\bfa_0-\bfb_0\|^2\right)}{\epsilon}}\right)$.
Then inequality~(\ref{eq:e}) and Jensen's inequality \cite{Rud87}
imply
\begin{equation}
  \label{eq:b}
  \forall t\geq T_\epsilon \quad
  \expect(\|\bfa(t)-\bfb(t)\|)\leq \sqrt{C/\lambda + \epsilon}    
\end{equation}    

Since $\|\bfa(t)-\bfb(t)\|$ is non-negative, (\ref{eq:b}) together
with Markov inequality \cite{Fel68} allow one to obtain the following
probabilistic bound on the distance between $\bfa(t)$ and $\bfb(t)$
\[
\forall A>0\ \forall t\geq T_\epsilon \quad  \proba\left(\|\bfa(t)-\bfb(t)\| \geq A \right) \leq
\frac{\sqrt{C/\lambda+\epsilon}}{A}
\]

Note however that this bound is much weaker than the asymptotic
almost-sure bound (\ref{eq:almost-sure}).

\subsubsection{Generalization to time-varying metrics}

\label{sec:gen}

Theorem \ref{theo:main} can be vastly generalized by considering
general time-dependent metrics (the case of state-dependent metrics is
not considered in this article and will be the subject of a future
work). Specifically, let us replace \textbf{(H1)} and \textbf{(H2)} by
the following hypotheses

\begin{description}
\item[(H1')] There exists a uniformly positive definite metric
  $\bfM(t)=\bfTh(t)^T\bfTh(t)$, with the lower-bound $\beta>0$ (i.e.
$
\forall \bfx,t \ \bfx^T\bfM(t)\bfx\geq \beta\|\bfx\|^2
$)
and $\bff(\bfa,t)$ is contracting in that metric, with
contraction rate $\lambda$, i.e.
\[
\lambda_{\max} \left(\left(\frac{d}{dt}\bfTh(t)+\bfTh(t)\frac{\partial
      \bff}{\partial \bfa}\right)\bfTh^{-1}(t)\right) \leq -\lambda
\quad \mathrm{uniformly}
\]
or equivalently
\[
\bfM(t)\frac{\partial \bff}{\partial \bfa}+\left(\frac{\partial \bff}{\partial
    \bfa}\right)^T 
\bfM(t)+\frac{d}{dt}\bfM(t)\leq -2\lambda \bfM(t) \quad \mathrm{uniformly}
\]
\item[(H2')] $\mathrm{tr}\left(\sigma(\bfa,t)^T\bfM(t)\sigma(\bfa,t)\right)$
  is uniformly upper-bounded by a constant $C$ 
\end{description}

\begin{defi}
  A system that verifies \emph{\textbf{(H1')}} and
  \emph{\textbf{(H2')}} is said to be \emph{stochastically
    contracting} in the metric $\bfM(t)$, with rate $\lambda$ and
  bound $C$.
\end{defi}

Consider now the generalized Lyapunov-like function
$V_1(\bfx,t)=(\bfa-\bfb)^T \bfM(t) (\bfa-\bfb)$.  Lemma
\ref{lemma:ineq} can then be generalized as follows.

\begin{lemma}
  \label{lemma:ineq2}
  Under \emph{\textbf{(H1')}} and \emph{\textbf{(H2')}}, one has the
  inequality
  \begin{equation}
    \widetilde{A}V_1(\bfx,t)\leq -2\lambda V_1(\bfx,t)+2C
  \end{equation}
\end{lemma}

\paragraph{Proof}

Let us compute first $\widetilde{A}V_1$

\begin{eqnarray}
\widetilde{A}V_1(\bfx,t)&=&\frac{\partial V_1}{\partial t}
  + \frac{\partial V_1}{\partial \bfx}\bffg(\bfx,t)
  +\frac{1}{2}\mathrm{tr}\left(
  \sg(\bfx,t)^T\frac{\partial^2 V_1}{\partial
      \bfx^2}\sg(\bfx,t)\right)    
   \nonumber\\ 
&=&(\bfa-\bfb)^T \left(\frac{d}{dt}\bfM(t)\right)(\bfa-\bfb) 
  +2(\bfa-\bfb)^T\bfM(t)(\bff(\bfa,t)-\bff(\bfb,t))\nonumber\\
&&+\mathrm{tr}(\sigma(\bfa,t)^T\bfM(t)\sigma(\bfa,t))
  +\mathrm{tr}(\sigma(\bfb,t)^T\bfM(t)\sigma(\bfb,t)) \nonumber 
\end{eqnarray}

Fix $t>0$ and consider the real-valued function
\[
r(\mu)=(\bfa-\bfb)^T\bfM(t)(\bff(\mu \bfa+(1-\mu)\bfb,t)-\bff(\bfb,t))
\]

Since $\bff$ is $C^1$, $r$ is $C^1$ over $[0,1]$. By the mean value
theorem, there exists $\mu_0\in]0,1[$ such that
\[
r'(\mu_0)=r(1)-r(0)=(\bfa-\bfb)^T\bfM(t)(\bff(\bfa)-\bff(\bfb))
\]

On the other hand, one obtains by differentiating $r$
\[
r'(\mu_0)=(\bfa-\bfb)^T\bfM(t)\left(\frac{\partial \bff}{\partial \bfa}(\mu_0
\bfa+(1-\mu_0)\bfb,t)\right)(\bfa-\bfb)
\]

Thus, letting $\bfc=\mu_0 \bfa+(1-\mu_0)\bfb$, one has
\begin{eqnarray}
\label{eq:bound-gen}
&(\bfa-\bfb)^T\left(\frac{d}{dt}\bfM(t)\right)(\bfa-\bfb)+2(\bfa-\bfb)^T\bfM(t)(\bff(\bfa)-\bff(\bfb))\nonumber \\
&=(\bfa-\bfb)^T\left(\frac{d}{dt}\bfM(t)\right)(\bfa-\bfb)+2(\bfa-\bfb)^T\bfM(t)\left(\frac{\partial
    \bff}{\partial \bfa}(\bfc,t)\right)(\bfa-\bfb) \nonumber \\
&=(\bfa-\bfb)^T \left(
  \frac{d}{dt}\bfM(t)+
  \bfM(t)\left(\frac{\partial \bff}{\partial \bfa}(\bfc,t)\right)+
  \left(\frac{\partial \bff}{\partial \bfa}(\bfc,t)\right)^T
  \bfM(t) \right) (\bfa-\bfb)\nonumber \\
&\leq -2\lambda (\bfa-\bfb)^T\bfM(t)(\bfa-\bfb)=-2\lambda V_1(\bfx)
\end{eqnarray}
where the inequality is obtained by using \textbf{(H1')}.

Finally, combining equation (\ref{eq:bound-gen}) with \textbf{(H2')}
allows to obtain the desired result. $\Box$

We can now state the generalized stochastic contraction theorem

\begin{theorem}[Generalized stochastic contraction]
  \label{theo:main-gen}
  Assume that system (\ref{eq:main}) verifies \emph{\textbf{(H1')}}
  and \emph{\textbf{(H2')}}. Let $\bfa(t)$ and $\bfb(t)$ be two
  trajectories whose initial conditions are given by a probability
  distribution $p(\bfx(0))=p(\bfa(0),\bfb(0))$. Then
  \begin{eqnarray}
    \label{eq:gen1}
    \forall t \geq 0 \quad \expect \left(
      (\bfa(t)-\bfb(t))^T\bfM(t)(\bfa(t)-\bfb(t)) \right) \leq
    \nonumber \\
    \frac{C}{\lambda} + e^{-2\lambda t}\int
      \left[(\bfa_0-\bfb_0)^T\bfM(0)(\bfa_0-\bfb_0)
        -\frac{C}{\lambda}\right]^+dp(\bfa_0,\bfb_0)
  \end{eqnarray}
  In particular,
  \begin{equation}
    \label{eq:gen2}
    \forall t \geq 0 \quad \expect \left(\|\bfa(t)-\bfb(t)\|^2 \right) \leq
    \frac{1}{\beta}\left(\frac{C}{\lambda} + \expect
      \left(\|\bfa(0)-\bfb(0)\|^2 \right)  e^{-2\lambda t}\right)
  \end{equation}

\end{theorem}

\paragraph{Proof} Following the same reasoning as in the proof of
theorem \ref{theo:main}, one obtains
\[
\forall t \geq 0 \quad \expectx V_1(\bfx(t)) \leq \frac{C}{\lambda}+
\left[V_1(\bfx_0)-\frac{C}{\lambda}\right]^+e^{-2\lambda t}
\]
which leads to (\ref{eq:gen1}) by integrating with respect to
$(\bfa_0,\bfb_0)$. Next, observing that
\[
\|\bfa(t)-\bfb(t)\|^2 \leq \frac{1}{\beta}(\bfa(t)-\bfb(t))^T\bfM(t)(\bfa(t)-\bfb(t))=
\frac{1}{\beta} \expect V_1(\bfx(t)) 
\]
and using the same bounding as in (\ref{eq:trick}) lead to
(\ref{eq:gen2}). $\Box$

\subsection{Strength of the stochastic contraction theorem}

\subsubsection{``Optimality'' of the mean square bound}

Consider the following linear dynamical system, known as the
Ornstein-Uhlenbeck (colored noise) process
\begin{equation}
  \label{eq:5}
  da=-\lambda a dt + \sigma dW
\end{equation}

Clearly, the noise-free system is contracting with rate $\lambda$ and
the trace of the noise matrix is upper-bounded by $\sigma^2$.  Let
$a(t)$ and $b(t)$ be two system trajectories starting respectively at
$a_0$ and $b_0$ (deterministic initial conditions). Then by theorem
\ref{theo:main}, we have
\begin{equation}
  \label{eq:optbound}
  \forall t \geq 0 \quad \expect\left((a(t)-b(t))^2 \right) \leq
  \frac{\sigma^2}{\lambda} +
  \left[(a_0-b_0)^2-\frac{\sigma^2}{\lambda}\right]^+e^{-2\lambda t}
\end{equation}

Let us verify this result by solving directly equation (\ref{eq:5}).
The solution of equation (\ref{eq:5}) is (\cite{Arn74}, p.~134)
\begin{equation}
  \label{eq:ousol}
  a(t)=a_0e^{-\lambda t}+\sigma\int_0^t e^{\lambda (s-t)}dW(s)
\end{equation}

Next, let us compute the mean square distance between the two
trajectories $a(t)$ and $b(t)$
\begin{eqnarray}
  \label{eq:1}
  \expect((a(t)-b(t))^2)&=&(a_0-b_0)^2e^{-2\lambda t} + \nonumber \\
  &&\sigma^2 \left( \expect\left(\left(\int_0^te^{\lambda (s-t)}dW_1(s)\right)^2
  \right)  + \expect\left(\left(\int_0^te^{\lambda (u-t)}dW_2(u)\right)^2
  \right) \right)\nonumber  \\
  &=&(a_0-b_0)^2e^{-2\lambda t} +
  \frac{\sigma^2}{\lambda}(1-e^{-2\lambda t}) \nonumber\\
  &\leq& \frac{\sigma^2}{\lambda} +
\left[(a_0-b_0)^2-\frac{\sigma^2}{\lambda}\right]^+e^{-2\lambda t} \nonumber
\end{eqnarray}

The last inequality is in fact an equality when $(a_0-b_0)^2\geq
\frac{\sigma^2}{\lambda}$.  Thus, this calculation shows that the
upper-bound (\ref{eq:optbound}) given by theorem \ref{theo:main} is
optimal, in the sense that it can be attained.

\subsubsection{No asymptotic almost-sure stability}
\label{sec:noalmostsure}

From the explicit form (\ref{eq:ousol}) of the solutions, one can
deduce that the distributions of $a(t)$ and $b(t)$ converge to the
normal distribution
$\mathscr{N}\left(0,\frac{\sigma^2}{2\lambda}\right)$
(\cite{Arn74},~p.~135). Since $a(t)$ and $b(t)$ are independent, the
distribution of the difference $a(t)-b(t)$ will then converge to
$\mathscr{N}\left(0,\frac{\sigma^2}{\lambda}\right)$. This observation
shows that, contrary to the case of standard stochastic stability (cf.
section~\ref{sec:standard-stab}), one cannot -- in general -- obtain
asymptotic almost-sure incremental stability results (which would
imply that the distribution of the difference converges instead to the
constant $0$).

Compare indeed equations (\ref{eq:ideal-bound}) (the condition for
standard stability, section~\ref{sec:standard-stab}) and
(\ref{eq:incremental-bound}) (the condition for incremental stability,
section~\ref{sec:stocon-theorem}). The difference lies in the term
$2C$, which stems from the fact that the influence of the noise does
not vanish when two trajectories get very close to each other (cf.
section~\ref{sec:settings}). The presence of this extra term prevents
$\widetilde{A}V(\bfx(t))$ from being always non-positive, and as a
result, it prevents $V(\bfx(t))$ from being always ``non-increasing''.
As a consequence, $V(\bfx(t))$ is not -- in general -- a
supermartingale, and one cannot then use the supermartingale
inequality to obtain asymptotic almost-sure bounds, as in equation
(\ref{eq:super}).

\paragraph{Remark} If one is interested in \emph{finite time} bounds
then the supermartingale inequality is still applicable, see
(\cite{Kus67}, p. 86) for details.

\subsection{Noisy and noise-free trajectories}
\label{sec:unperturbed}
 
Consider the following augmented system
\begin{equation}
  \label{eq:dup2}
  d\bfx=
  \left(\begin{array}{l}
      \bff(\bfa,t)\\
      \bff(\bfb,t)
    \end{array}\right)dt+
  \left(\begin{array}{cc}
      0&0\\
      0&\sigma(\bfb,t)
    \end{array}\right)
  \left(\begin{array}{c}
      dW_d^1\\
      dW_d^2
    \end{array}\right)=\bffg(\bfx,t)dt+\sg(\bfx,t)dW_{2d}
\end{equation}

This equation is the same as equation (\ref{eq:dup}) except that the
$\bfa$-system is not perturbed by noise. Thus
$V(\bfx)=\|\bfa-\bfb\|^2$ will represent the distance between a
noise-free trajectory and a noisy one. All the calculations will be
the same as in the previous development, with $C$ being replaced by
$C/2$.  One can then derive the following corollary

\begin{cor}
  \label{cor:unpert}
  Assume that system (\ref{eq:main}) verifies \emph{\textbf{(H1')}}
  and \emph{\textbf{(H2')}}. Let $\bfa(t)$ be a \emph{noise-free}
  trajectory starting at $\bfa_0$ and $\bfb(t)$ a \emph{noisy}
  trajectory whose initial condition is given by a probability
  distribution $p(\bfb(0))$. Then
  \begin{equation}
    \label{eq:e2} 
    \forall t \geq 0  \quad \expect\left(\|\bfa(t)-\bfb(t)\|^2 \right) \leq
    \frac{1}{\beta}\left(\frac{C}{2\lambda}+
    \expect\left(\|\bfa_0-\bfb(0)\|^2\right)e^{-2\lambda t}\right)
  \end{equation}
\end{cor}

\paragraph{Remarks} 
\begin{itemize}
\item One can note here that the derivation of corollary
  \ref{cor:unpert} is only permitted by our initial choice of
  considering \emph{distinct} driving Wiener process for the $\bfa$-
  and $\bfb$-systems (cf.  section~\ref{sec:settings}).
\item Corollary \ref{cor:unpert} provides a \emph{robustness} result
  for contracting systems, in the sense that any contracting system is
  \emph{automatically} protected against noise, as quantified by
  (\ref{eq:e2}). This robustness could be related to the exponential
  nature of contraction stability.
\end{itemize}


%% file: combinations.tex
\section{Combinations of contracting stochastic systems}
\label{sec:combinations}

Stochastic contraction inherits naturally from deterministic
contraction \cite{LS98} its convenient combination properties.
Because contraction is a state-space concept, such properties can be
expressed in more general forms than input-output analogues such as
passivity-based combinations~\cite{Pop73}. The following combination
properties allow one to build by recursion stochastically contracting
systems of arbitrary size.

\paragraph{Parallel combination}

Consider two stochastic systems of the same dimension
\[
\left\{
\begin{array}{l}
d\bfx_1=\bff_1(\bfx_1,t)dt+\sigma_1(\bfx_1,t)dW_1\\
d\bfx_2=\bff_2(\bfx_2,t)dt+\sigma_2(\bfx_2,t)dW_2
\end{array}\right.
\]

Assume that both systems are stochastically contracting in the same
\emph{constant} metric $\bfM$, with rates $\lambda_1$ and $\lambda_2$
and with bounds $C_1$ and $C_2$. Consider a uniformly positive bounded
superposition
\[
\alpha_1(t)\bfx_1+\alpha_2(t)\bfx_2
\]
where $\forall t\geq 0,\
l_i\leq\alpha_i(t)\leq m_i$ for some $l_i,m_i>0,\ i=1,2$.

Clearly, this superposition is stochastically contracting in the
metric $\bfM$, with rate $l_1\lambda_1+l_2\lambda_2$ and bound
\mbox{$m_1C_1+m_2C_2$}.

\paragraph{Negative feedback combination}

In this and the following paragraphs, we describe combinations
properties for contracting systems in constant metrics $\bfM$. The
case of time-varying metrics can be easily adapted from this
development but is skipped here for the sake of clarity.

Consider two coupled stochastic systems
\[
\left\{
\begin{array}{l}
d\bfx_1=\bff_1(\bfx_1,\bfx_2,t)dt+\sigma_1(\bfx_1,t)dW_1\\
d\bfx_2=\bff_2(\bfx_1,\bfx_2,t)dt+\sigma_2(\bfx_2,t)dW_2
\end{array}\right.
\]
Assume that system $i$ ($i=1,2$) is stochastically contracting with
respect to $\bfM_i=\bfTh_i^T\bfTh_i$, with rate $\lambda_i$ and bound
$C_i$.

Assume furthermore that the two systems are connected by
\emph{negative feedback} \cite{TS05}. More precisely, the Jacobian
matrices of the couplings are of the form
$\bfTh_1\bfJ_{12}\bfTh_2^{-1}=-k\bfTh_2\bfJ_{21}^T\bfTh_1^{-1}$, with
$k$ a positive constant. Hence, the Jacobian matrix of the augmented
system is given by
\[\bfJ=\left(
\begin{array}{cc}
\bfJ_1&-k\bfTh_1^{-1}\bfTh_2\bfJ_{21}^T\bfTh_1^{-1}\bfTh_2\\
\bfJ_{21}&\bfJ_2
\end{array}
\right)
\]

Consider a coordinate transform
$
\bfTh=\left(
\begin{array}{cc}
\bfTh_1&\zeros\\
\zeros&\sqrt{k}\bfTh_2
\end{array}
\right)
$
associated with the metric $\bfM=\bfTh^T\bfTh>\zeros$. After some
calculations, one has

\begin{eqnarray}
\left(\bfTh\bfJ\bfTh^{-1}\right)_s&=&
\left(
\begin{array}{cc}
\left(\bfTh_1\bfJ_1\bfTh_1^{-1}\right)_s&\zeros\\
\zeros&\left(\bfTh_2\bfJ_2\bfTh_2^{-1}\right)_s
\end{array}
\right)\nonumber\\
&\leq&\max(-\lambda_1,-\lambda_2)\bfI\quad
{\rm uniformly}
\end{eqnarray}

The augmented system is thus stochastically contracting in the metric
$\bfM$, with rate $\min(\lambda_1,\lambda_2)$ and bound $C_1+kC_2$.

\paragraph{Hierarchical combination}

We first recall a standard result in matrix analysis
\cite{HJ85}. Let $\bfA$ be a symmetric matrix in the form $
\bfA= \left(\begin{array}{cc}
    \bfA_1&\bfA_{21}^T\\
    \bfA_{21}&\bfA_2
\end{array}\right)
$.
Assume that $\bfA_1$ and $\bfA_2$ are definite positive. Then $\bfA$
is definite positive if
$
\mathrm{sing}^2(\bfA_{21})<\lambda_{\min}(\bfA_1)\lambda_{\min}
(\bfA_2)
$
where $\mathrm{sing}(\bfA_{21})$ denotes the largest singular value
of $\bfA_{21}$. In this case, the smallest eigenvalue of $\bfA$
satisfies
\[
\lambda_{\min}(\bfA)\geq
\frac{\lambda_{\min}(\bfA_1)+\lambda_{\min}(\bfA_2)}{2}
-\sqrt{\left(\frac{\lambda_{\min}(\bfA_1)-\lambda_{\min}(\bfA_2)}{2}\right)^2+
\mathrm{sing}^2(\bfA_{21})}
\]

Consider now the same set-up as in the previous paragraph, except that
the connection is now \emph{hierarchical} and upper-bounded. More
precisely, the Jacobians of the couplings verify $\bfJ_{12}=\zeros$
and $\mathrm{sing}^2(\bfTh_2\bfJ_{21}\bfTh_1^{-1})\leq K$.  The
Jacobian matrix of the augmented system is then given by 
$ \bfJ=
\left(\begin{array}{cc}
    \bfJ_1&\zeros\\
    \bfJ_{21}&\bfJ_2
\end{array}\right)
$.
Consider a coordinate transform
$
\bfTh_\epsilon=\left(
\begin{array}{cc}
\bfTh_1&\zeros\\
\zeros&\epsilon\bfTh_2
\end{array}
\right)
$ 
associated with the metric
$\bfM_\epsilon=\bfTh_\epsilon^T\bfTh_\epsilon>\zeros$. After some
calculations, one has
\begin{eqnarray}
\left(\bfTh\bfJ\bfTh^{-1}\right)_s&=&
\left(
\begin{array}{cc}
\left(\bfTh_1\bfJ_1\bfTh_1^{-1}\right)_s&
\frac{1}{2}\epsilon(\bfTh_2\bfJ_{21}\bfTh_1^{-1})^T\\
\frac{1}{2}\epsilon\bfTh_2\bfJ_{21}\bfTh_1^{-1}&
\left(\bfTh_2\bfJ_2\bfTh_2^{-1}\right)_s
\end{array}
\right)\nonumber
\end{eqnarray}
Set now $\epsilon=\sqrt{\frac{2\lambda_1\lambda_2}{K}}$. The augmented
system is then stochastically contracting in the metric
$\bfM_\epsilon$, with rate
$\frac{1}{2}(\lambda_1+\lambda_2-\sqrt{\lambda_1^2+\lambda_2^2}))$ and
bound $C_1+\frac{2C_2\lambda_1\lambda_2}{K}$.

\paragraph{Small gains}

In this paragraph, we require no specific assumption on the form of
the couplings. Consider the coordinate transform 
$ \bfTh=\left(
\begin{array}{cc}
\bfTh_1&\zeros\\
\zeros&\sqrt{k}\bfTh_2
\end{array}
\right)
$
associated with the metric $\bfM_k=\bfTh_k^T\bfTh_k>\zeros$. Aftersome
calculations, one has

\begin{eqnarray}
\left(\bfTh_k\bfJ\bfTh_k^{-1}\right)_s&=&
\left(
\begin{array}{cc}
\left(\bfTh_1\bfJ_1\bfTh_1^{-1}\right)_s&
\bfB_k^T\\
\bfB_k&
\left(\bfTh_2\bfJ_2\bfTh_2^{-1}\right)_s
\end{array}
\right)\nonumber
\end{eqnarray}
where $\bfB_k=\frac{1}{2}\left(\sqrt{k}\bfTh_2\bfJ_{21}\bfTh_1^{-1}+
\frac{1}{\sqrt{k}}\left(\bfTh_1\bfJ_{12}\bfTh_2^{-1}\right)^T\right)$.

Following the matrix analysis result stated at the beginning of the
previous paragraph, if $\inf_{k>0}\mathrm{sing}^2(\bfB_k)<\lambda_1
\lambda_2$ then the augmented system is stochastically contracting in
the metric $\bfM_k$, with bound $C_1+kC_2$ and rate $\lambda$ verifying
\begin{equation}
  \lambda\geq\frac{\lambda_1+\lambda_2}{2}
  -\sqrt{\left(\frac{\lambda_1-\lambda_2}{2}\right)^2+
    \inf_{k>0}\mathrm{sing}^2(\bfB_k)}
\end{equation}


%% file: examples.tex
\section{Some examples}
\label{sec:examples}

\subsection{Effect of measurement noise on contracting observers}

Consider a nonlinear dynamical system
\begin{equation}
  \label{eq:toto}
  \dot\bfx = \bff(\bfx,t)  
\end{equation}

If a measurement $\bfy= \bfy (\bfx)$ is available, then it may be
possible to choose an output injection matrix $\bfK(t)$ such that the
dynamics
\begin{equation}
  \label{eq:det-obs}
\dot{\hat\bfx} = \bff(\hat\bfx,t) + \bfK(t)(\hat\bfy - \bfy)  
\end{equation}
is contracting, with $\hat\bfy=\bfy(\hat\bfx)$. Since the actual state
$\bfx$ is a particular solution of (\ref{eq:det-obs}), any solution
$\hat\bfx$ of (\ref{eq:det-obs}) will then converge towards $\bfx$
exponentially.

Assume now that the measurements are corrupted by additive ``white
noise''. In the case of \emph{linear} measurement, the measurement
equation becomes \mbox{$\bfy=\bfH(t)\bfx+\Sigma(t)\xi(t)$} where
$\xi(t)$ is a multidimensional ``white noise'' and $\Sigma(t)$ is the
matrix of measurement noise intensities.

The observer equation is now given by the following It\^o stochastic
differential equation (using the formal rule $dW=\xi dt$)
\begin{equation}
  \label{eq:stochastic-observer}
  d{\hat{\bfx}}=(\bff(\hat{\bfx},t)+\bfK(t)(\bfH(t)\bfx-\bfH(t)\hat{\bfx}))dt
  + \bfK(t)\Sigma(t) dW
\end{equation}

Next, remark that the solution $\bfx$ of system (\ref{eq:toto}) is a
also a solution of the noise-free version of system
(\ref{eq:stochastic-observer}). By corollary~\ref{cor:unpert}, one
then has, for any solution $\hat\bfx$ of system
(\ref{eq:stochastic-observer})
\begin{equation}
  \label{eq:obs-bound}
  \forall t \geq 0 \quad \expect\left(\|\hat{\bfx}(t)-\bfx(t)\|^2
  \right) \leq \frac{C}{2\lambda}+
  \|\hat{\bfx}_0-\bfx_0\|^2e^{-2\lambda
    t}
\end{equation}
where 
\[
\lambda= \inf_{\bfx,t}\left|\lambda_{\max}\left(
    \frac{\partial{\bff(\bfx,t)}}{\partial\bfx}-\bfK(t)\bfH(t) \right) \right|
\]
\[
C=\sup_{t\geq 0} \mathrm{tr}\left(\Sigma(t)^T\bfK(t)^T\bfK(t)\Sigma(t)\right)
\]

\paragraph{Remark} The choice of the injection gain $\bfK(t)$ is
governed by a trade-off between convergence speed ($\lambda$) and
noise sensitivity ($C/\lambda$) as quantified by~(\ref{eq:obs-bound}).
More generally, the explicit computation of the bound on the expected
quadratic estimation error given by (\ref{eq:obs-bound}) may open the
possibility of {\it measurement selection} in a way similar to the
linear case. If several possible measurements or sets of measurements
can be performed, one may try at each instant (or at each step, in a
discrete version) to select the most relevant, i.e., the measurement
or set of measurements which will best contribute to improving the
state estimate. Similarly to the Kalman filters used in~\cite{Dic98}
for linear systems, this can be achieved by computing, along with the
state estimate itself, the corresponding bounds on the expected
quadratic estimation error, and then selecting accordingly the
measurement which will minimize it.

\subsection{Estimation of velocity using composite variables}

In this section, we present a very simple example that hopefully
suggests the many possibilities that could stem from the combination
of our stochastic stability analysis with the composite variables
framework \cite{SL91}.

Let $x$ be the position of a mobile subject to a sinusoidal forcing
\[
\ddot{x}=-U_1\omega^2\sin(\omega t)+2U_2
\]
where $U_1$ and $\omega$ are \emph{known} parameters. We would like to
compute good approximations of the mobile's velocity $v$ and
acceleration $a$ using only measurements of $x$ and without using any
filter. For this, construct the following observer
\begin{eqnarray}
  \label{eq:composite-obs}
  \frac{d}{dt}
\left(
\begin{array}{c}
  \overline{v} \\
  \overline{a}
\end{array}
\right)
&=&
\left(
\begin{array}{cc}
  -\alpha_v & 1 \\
  -\alpha_a & 0 
\end{array}
\right)
\left(
\begin{array}{c}
  \overline{v} \\
  \overline{a}
\end{array}
\right)
+
\left(
\begin{array}{c}
  (\alpha_a-\alpha_v^2)x \\
  -\alpha_a \alpha_v x - U_1\omega^3\cos(\omega t)
\end{array}
\right) \nonumber \\
&=&
\bfA
\left(
\begin{array}{c}
  \overline{v} \\
  \overline{a}
\end{array}
\right)
+
\bfB
\left(
\begin{array}{c}
  x \\
  x
\end{array}
\right)
-\left(
\begin{array}{c}
  0 \\
  U_1\omega^3\cos(\omega t)
\end{array}
\right)
\end{eqnarray}
and introduce the composite variables $\widehat{v} =
\overline{v} + \alpha_v{x}$ and $\widehat{a} = \overline{a}
+ \alpha_a{x}$. By construction, these variables follow the equation 
\begin{equation}
  \label{eq:composite}
  \frac{d}{dt}
\left(
\begin{array}{c}
  \widehat{v} \\
  \widehat{a}
\end{array}
\right)
=
\bfA
\left(
\begin{array}{c}
  \widehat{v} - v \\
  \widehat{a}
\end{array}
\right)
+
\left(
\begin{array}{c}
  0 \\
  - U_1\omega^3\cos(\omega t)
\end{array}
\right)
\end{equation}
and therefore, a particular solution of $(\widehat{v},\widehat{a})$ is
clearly $(v,a)$. Choose now $\alpha_a = \alpha_v^2=\alpha^2$ and let
$\bfM_\alpha=\frac{1}{2}\left(
  \begin{array}{cc}
    \alpha^2&-\alpha/2\\
    -\alpha/2&1
  \end{array}
\right)$. One can then show that system (\ref{eq:composite}) is
contracting with rate $\lambda_\alpha=\alpha/2$ in the metric
$\bfM_\alpha$.  Thus, by the basic contraction theorem \cite{LS98},
$(\widehat{v},\widehat{a})$ converges exponentially to $(v,a)$ with
rate $\lambda_\alpha$ in the metric $\bfM_\alpha$. Also note that the
$\beta$-bound corresponding to the metric $\bfM_\alpha$ is given by
$\beta_\alpha=\frac{1+\alpha^2-\sqrt{\alpha^4-\alpha^2+1}}{4}$.

Next, assume that the measurements of $x$ are corrupted by additive
``white noise'', so that $x_\mathrm{measured}=x+\sigma\xi$.  Equation
(\ref{eq:composite-obs}) then becomes an It\^o stochastic differential
equation
\[
d
\left(
\begin{array}{c}
  \overline{v} \\
  \overline{a}
\end{array}
\right)
=
\left[
\bfA
\left(
\begin{array}{c}
  \overline{v} \\
  \overline{a}
\end{array}
\right)
+
\bfB
\left(
\begin{array}{c}
  x \\
  x
\end{array}
\right)
-\left(
\begin{array}{c}
  0 \\
  U_1\omega^3\cos(\omega t)
\end{array}
\right)
\right]
dt
+
\bfB \left(
\begin{array}{cc}
  \sigma&0 \\
  0 &\sigma
\end{array}
\right) dW
\]

\begin{figure}[h!] 
  \begin{minipage}[ht]{0.5\linewidth}
    \centering
    \includegraphics[scale=0.5]{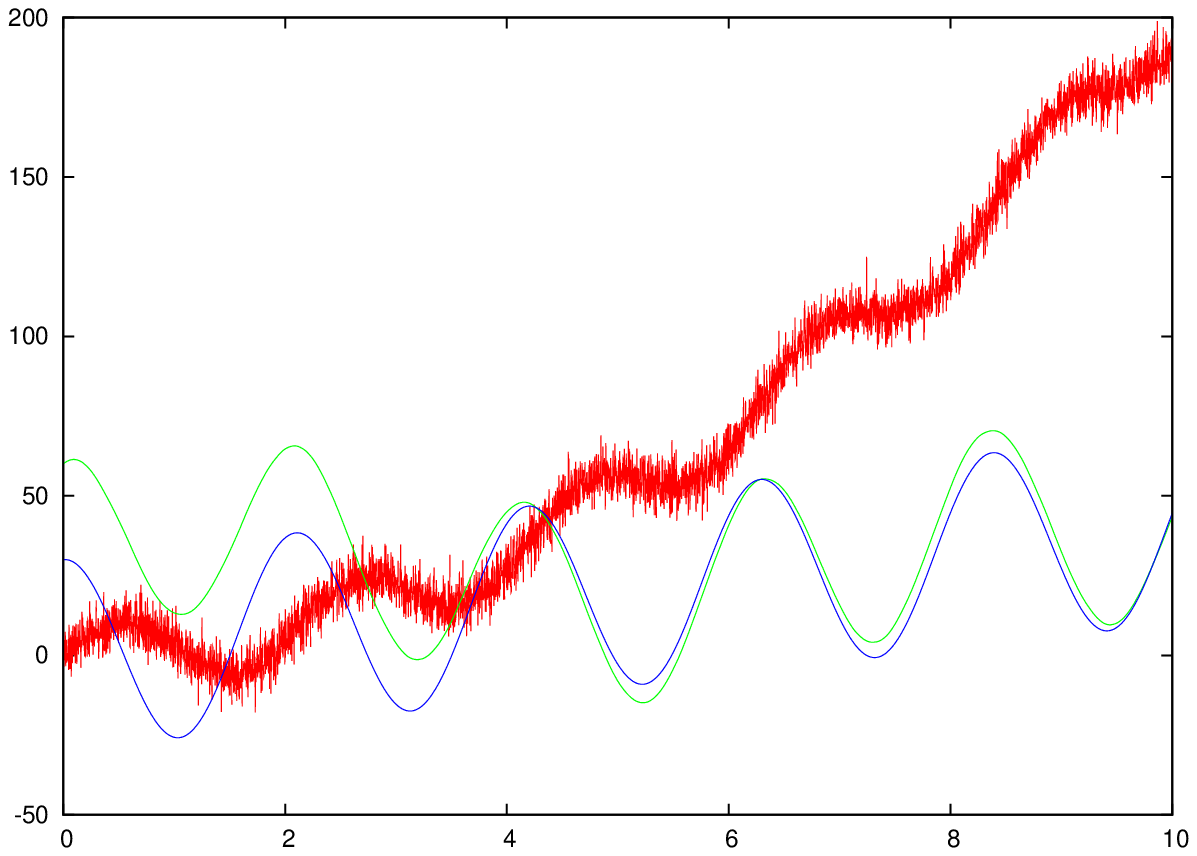}
    \end{minipage}
  \begin{minipage}[ht]{0.5\linewidth}
    \centering
    \includegraphics[scale=0.5]{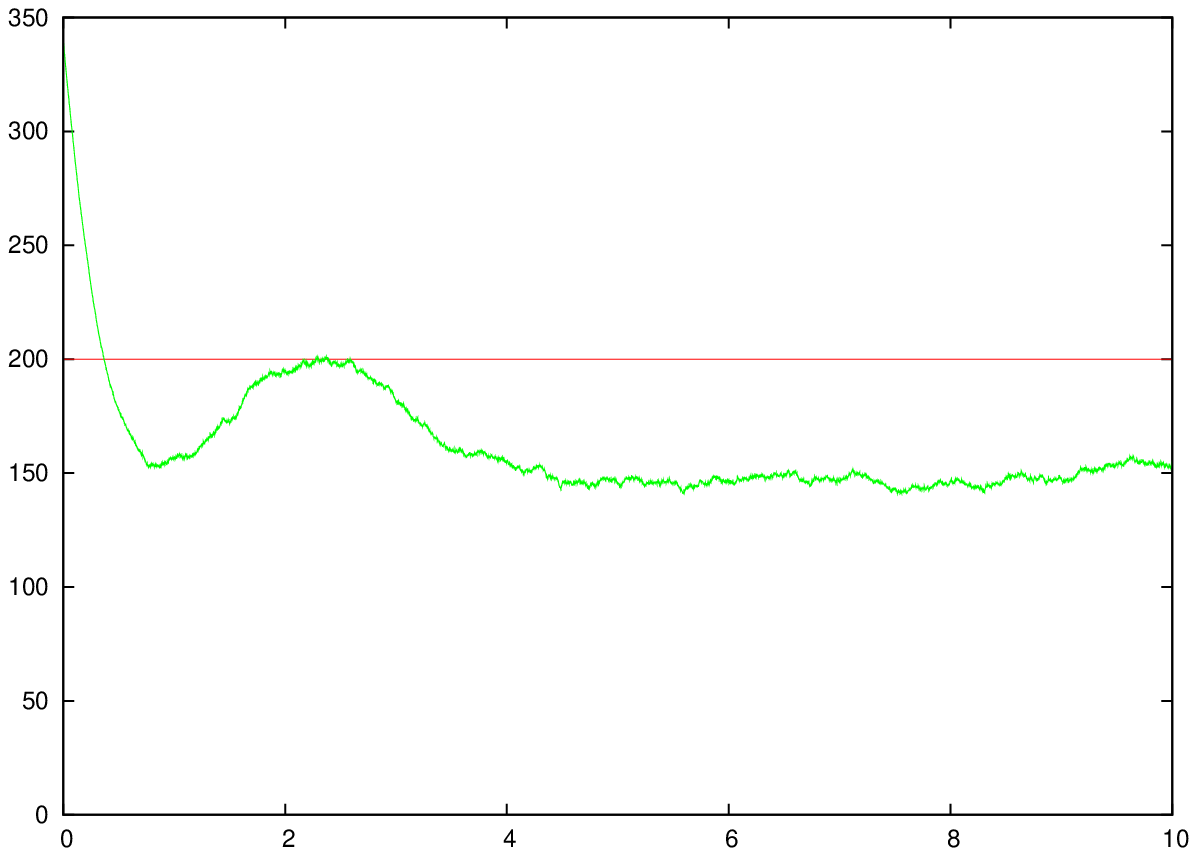}
  \end{minipage}
  \caption[]{Estimation of the velocity of a mobile using noisy
    measurements of its position. The simulation was performed using
    the Euler-Maruyama algorithm \cite{Hig01} with the following
    parameters: $U_1=10$, $U_2=2$, $\omega = 3$, $\sigma=10$ and
    $\alpha =1$.  Left plot: simulation for one trial. The plot shows
    the measured position (red), the actual velocity (blue) and the
    estimate of the velocity using the measured position (green).
    Right plot: the average over 1000 trials of the squared error
    $\|\widehat{v}-v\|^2+\|\widehat{a}-a\|^2$ (green) and the
    asymptotic bound
    $\left(\frac{\alpha^5\sigma^2}{2\beta_\alpha}=200\right)$ given by
    our approach (red).}
  \label{fig:composite_variable}
\end{figure}

By definition of $\bfB$, the variance of the noise in the metric
$\bfM_\alpha$ is upper-bounded by $\frac{\alpha^6\sigma^2}{2}$. Thus,
using again corollary \ref{cor:unpert}, one obtains (see Figure
\ref{fig:composite_variable} for a numerical simulation)
\[ 
\forall t \geq 0 \quad \expect\left(\|\widehat{v}(t)-v(t)\|^2 +
  \|\widehat{a}(t)-a(t)\|^2 \right) \leq \frac{\alpha^5\sigma^2}{2\beta_\alpha} +
   \frac{\|\widehat{v}_0-v_0\|^2+\|\widehat{a}_0-a_0\|^2}{2\beta_\alpha}
  e^{-\alpha t}
\]

\subsection{Stochastic synchronization}
\label{sec:sync}

Consider a network of $n$ dynamical elements coupled through diffusive
connections
\begin{equation}
  \label{equ:robustsync}
  d\bfx_i = \left(\bff(\bfx_i,t) + \sum_{j \neq i} \bfK_{ij}(\bfx_j -
  \bfx_i)\right)dt + \sigma_i(\bfx_i,t)dW^d_i \qquad i=1,\ldots,n
\end{equation}

Let
\[
\bfxg=
\left(\begin{array}{c}
    \bfx_1\\
    \vdots\\
    \bfx_n
\end{array}\right), \quad
\bffg(\bfxg,t)=
\left(\begin{array}{c}
    \bff(\bfx_1,t)\\
    \vdots\\
    \bff(\bfx_n,t)
\end{array}\right), \quad
\sg(\bfxg,t)=
\left(\begin{array}{ccc}
\sigma_1(\bfx_1,t)&\zeros&\zeros\\
\zeros&\ddots&\zeros\\
\zeros&\zeros&\sigma_n(\bfx_n,t)
\end{array}\right)
\]

The global state $\bfxg$ then follows the equation
\begin{equation}
  \label{eq:sync-global}
  d\bfxg=\left(\bffg(\bfxg,t)-\bfL\bfxg\right)dt+\sg(\bfxg,t)dW^{nd}  
\end{equation}

In the sequel, we follow the reasoning of \cite{PS07}, which starts
by defining an appropriate orthonormal matrix $\bfV$ describing the
synchronization subspace ($\bfV$ represents the state projection on
the subspace $\sM^\perp$, orthogonal to the synchronization subspace
$\sM=\{(\bfx_1,\dots,\bfx_n)^T:\bfx_1=\dots=\bfx_n\}$, see
\cite{PS07} for details).  Denote by $\bfyg$ the state of the
projected system, $\bfyg=\bfV\bfxg$. Since the mapping is linear,
It\^o differentiation rule simply yields
\begin{eqnarray}
  \label{eq:sync}
  d\bfyg&=&\bfV d\bfxg= \left(\bfV\bffg(\bfxg,t)-\bfV\bfL\bfxg\right)dt
  +\bfV\sg(\bfxg,t)dW^{nd} \nonumber\\
  &=&\left(\bfV\bffg(\bfV^T\bfyg,t)-\bfV\bfL\bfV^T\bfyg\right)dt
  +\bfV\sg(\bfV^T\bfyg,t)dW^{nd}
\end{eqnarray}

Assume now that $\frac{\partial{\bff}}{\partial{\bfx}}$ is uniformly
upper-bounded. Then for strong enough coupling strength,
$\bfA=\bfV\frac{\partial{\bff}}{\partial{\bfx}}\bfV^T-\bfV\bfL\bfV^T$
will be uniformly negative definite. Let $\lambda =
|\lambda_\mathrm{max}(\bfA)| >0$. System (\ref{eq:sync}) then verifies
condition \textbf{(H1)} with rate $\lambda$. Assume furthermore that
each noise intensity $\sigma_i$ is upper-bounded by a constant $C_i$
(i.e.  $\sup_{\bfx,t}
\mathrm{tr}(\sigma_i(\bfx,t)^T\sigma_i(\bfx,t))\leq C_i$). Condition
\textbf{(H2)} will then be satisfied with the bound $C=\sum_{i}C_i$.

Next, consider a noise-free trajectory $\bfyg_u(t)$ of system
(\ref{eq:sync}). By theorem~3 of \cite{PS07}, we know that
$\bfyg_u(t)$ converges exponentially to zero. Thus, by corollary
\ref{cor:unpert}, one can conclude that, after exponential transients
of rate $\lambda$, $ \expect(\|\bfyg(t)\|^2)\leq \frac{C}{2\lambda}
$.

On the other hand, one can show that
\[
\|\bfyg(t)\|^2=\frac{1}{n}\sum_{i,j}\|\bfx_i-\bfx_j\|^2
\]

Thus, after exponential transients of rate $\lambda$, we have
\[
\sum_{i,j}\|\bfx_i-\bfx_j\|^2
\leq  \frac{nC}{2\lambda}
\]

\paragraph{Remarks}
\begin{itemize}
\item The above development is fully compatible with the concurrent
  synchronization framework \cite{PS07}. It can also be easily
  generalized to the case of time-varying metrics by combining
  theorem \ref{theo:main-gen} of this paper and corollary~1 of
  \cite{PS07}.
\item The synchronization of It\^o dynamical systems has been
  investigated in \cite{CK05sync}. However, the systems considered by
  the authors of that article were \emph{dissipative}. Here, we make a
  less restrictive assumption, namely, we only require
  $\frac{\partial{\bff}}{\partial{\bfx}}$ to be uniformly
  upper-bounded. This enables us to study the synchronization of a
  broader class of dynamical systems, which can include nonlinear
  oscillators or even chaotic systems.
\end{itemize}

\paragraph{Example} As illustration of the above development, we
provide here a detailed analysis for the synchronization of noisy
FitzHugh-Nagumo oscillators (see \cite{WS05} for the references).
The dynamics of two diffusively-coupled noisy FitzHugh-Nagumo
oscillators can be described by
\[
\left\{\begin{array}{l}
    dv_i=(c(v_i+w_i-\frac{1}{3} v_i^3 +I_i)+k(v_0-v_i))dt +\sigma dW_i \\
    dw_i=-\frac{1}{c}(v_i-a +b w_i)dt
\end{array}\right. \quad i=1,2
\]

Let $\bfx=(v_1,w_1,v_2,w_2)^T$ and 
$
\bfV=\frac{1}{\sqrt 2}
\left(
\begin{array}{cccc}
  1&0&-1&0\\0&1&0&-1
\end{array}\right)
$. 
The Jacobian matrix of the projected noise-free system is then given
by
\[
\left(
\begin{array}{cc}
  c-\frac{c(v_1^2+v_2^2)}{2}-k& c \\
  -1/c&-b/c
\end{array}\right)
\]

\begin{figure}[h!] 
  \begin{center}
    \includegraphics[scale=0.5]{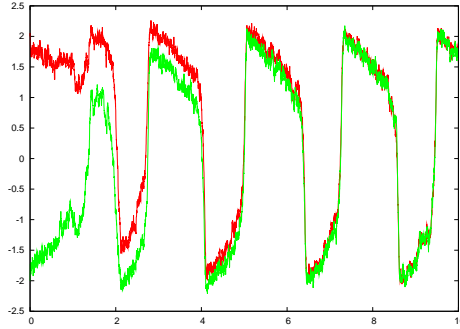}
    \caption{Synchronization of two noisy FitzHugh-Nagumo oscillators.
      The simulation was performed using the Euler-Maruyama algorithm
      \cite{Hig01} with the following parameters: $a=0.3$, $b=0.2$,
      $c=30$, $k=40$ and $\sigma=1$. The plot shows the ``membrane
      potentials'' of the two oscillators.}
    \label{simu:fn}
  \end{center}
\end{figure}

Thus, if the coupling strength verifies $k>c$ then the projected
system will be stochastically contracting in the diagonal metric
$\bfM=\mathrm{diag}(1,c)$ with rate $\min(k-c,b/c)$ and bound
$\sigma^2$. Hence, the average absolute difference between the two
``membrane potentials'' $|v_1-v_2|$ will be upper-bounded by
${\sigma}/{\sqrt{{\min(1,c)\min(k-c,b/c)}}}$ (see Figure \ref{simu:fn}
for a numerical simulation).

%% file: appendix.tex
\appendix
\section{Appendix}

\subsection{Proof of the supermartingale property}

\begin{lemma}
  \label{lemma:super}
  Consider a Markov stochastic process $\bfx(t)$ and a non-negative
  function $V$ such that \mbox{$\forall t\geq 0 \quad \expect
    V(\bfx(t)) < \infty$} and
  \begin{equation}
    \forall \bfx\in\mathbb{R}^n \quad \widetilde{A}V(\bfx) \leq -\lambda
    V(\bfx)
  \end{equation}  
  where $\lambda$ is a non-negative real number and $\widetilde{A}$ is
  the infinitesimal operator of the process $\bfx(t)$. Then
  $V(\bfx(t))$ is a supermartingale with respect to the canonical
  filtration $\mathscr{F}_t=\{\bfx(s),s\leq t\}$.
\end{lemma}

We need to show that for all $s\geq t$, one has
$\expect(V(\bfx(s))|\mathscr{F}_t) \leq V(\bfx(t))$. Since $\bfx(t)$ is
a Markov process, it suffices to show that
\[
\forall \bfx_0\in \mathbb{R}^n \quad
\expect(V(\bfx(t))|\bfx(0)=\bfx_0)\leq V(\bfx_0)
\]

By Dynkin's formula, one has for all $\bfx_0\in \mathbb{R}^n$
\[
\begin{array}{rcl}
  \expect_{\bfx_0}V(\bfx(t))&=&V(\bfx_0)+\expect_{\bfx_0}\int_0^t
  \widetilde{A}V(\bfx(s)) ds\\
  &\leq& V(\bfx_0)-\lambda \expect_{\bfx_0} \int_0^t V(\bfx(s))ds \leq V(\bfx_0)
\end{array}
\]
where $\expect_{\bfx_0}(\cdot)=\expect(\cdot|\bfx(0)=\bfx_0)$.

\subsection{A variation of Gronwall's lemma}

\begin{lemma}
  \label{lemma:gronwall}
  Let $g:[0,\infty[ \to\mathbb{R}$ be a continuous function, $C$ a
  real number and $\lambda$ a \emph{strictly positive} real
  number. Assume that
  \begin{equation}
    \label{eq:bound1}
    \forall u,t \quad 0\leq u \leq t  \quad  g(t)-g(u) \leq \int_u^t -\lambda
    g(s)+C ds
  \end{equation}
  Then
  \begin{equation}
    \label{eq:bound2}
    \forall t \geq 0 \quad g(t)\leq \frac{C}{\lambda}+
    \left[g(0)-\frac{C}{\lambda}\right]^+e^{-\lambda t}
  \end{equation}
  where $[\cdot]^+=\max(0,\cdot)$.
\end{lemma}

\paragraph{Proof}

\textbf{Case 1 :} $C=0$, $g(0)>0$.

Define $h(t)$ by
\[
\forall t \geq 0 \quad  h(t)=g(0)e^{-\lambda t}
\]
Remark that $h$ is positive with $h(0)=g(0)$, and satisfies
(\ref{eq:bound1}) where the inequality has been replaced by an
equality
\[
\forall u,t \quad 0\leq u \leq t  \quad h(t)-h(u) = -\int_u^t \lambda
h(s) ds
\]

Consider now the set 
$
S=\{t\geq 0 \ |\  g(t) > h(t) \}
$.
If $S=\emptyset$ then the lemma holds true. Assume by contradiction
that $S\neq \emptyset$. In this case, let $m=\inf S <\infty$. By
continuity of $g$ and $h$ and by the fact that $g(0)=h(0)$, one has
$g(m)=h(m)$ and there exists $\epsilon>0$ such that
\begin{equation}
  \label{eq:epsilon}
  \forall t\in ]m,m+\epsilon[ \quad g(t)>h(t)  
\end{equation}

Consider now $\phi(t)=g(m)-\lambda\int_m^t g(s)ds$. Equation 
(\ref{eq:bound1}) implies that
\[
\forall t\geq m \quad g(t)\leq \phi(t)
\]
In order to compare $\phi(t)$ and $h(t)$ for
$t\in ]m,m+\epsilon[$, let us differentiate the ratio $\phi(t)/h(t)$.
\[
\left(\frac{\phi}{h}\right)'=\frac{\phi'h-h'\phi}{h^2}
=\frac{-\lambda gh+\lambda h \phi}{h^2}=\frac{\lambda
  h(\phi-g)}{h^2}\geq 0
\]
Thus $\phi(t)/h(t)$ is increasing for $t\in ]m,m+\epsilon[$. Since
$\phi(m)/h(m)=1$, one can conclude that
\[
\forall t\in ]m,m+\epsilon[ \quad \phi(t)\geq h(t)
\]
which implies, by definition of $\phi$ and $h$, that
\begin{equation}
  \label{eq:t}
  \forall t\in ]m,m+\epsilon[ \quad \int_m^t g(s)ds \leq \int_m^t h(s)ds  
\end{equation}

Choose now $t_0$ such that $m<t_0<m+\epsilon$, then one has by 
(\ref{eq:epsilon}) 
\[
\int_m^{t_0}g(s)ds > \int_m^{t_0}h(s)ds
\]
which clearly contradicts (\ref{eq:t}).

\textbf{Case 2 :} $C=0$, $g(0)\leq0$

Consider the set 
$
S=\{t\geq 0 \ | \  g(t) > 0 \}
$.
If $S=\emptyset$ then the lemma holds true. Assume by contradiction
that $S\neq \emptyset$. In this case, let $m=\inf S<\infty$. By
continuity of $g$ and by the fact that $g(0)\leq 0$, one has $g(m)=0$
and there exists $\epsilon$ such that
\begin{equation}
  \label{eq:epsilon2}
  \forall t\in ]m,m+\epsilon[ \quad g(t)>0  
\end{equation}
Let $t_0\in ]m,m+\epsilon[$. Equation (\ref{eq:bound1}) implies that
\[
g(t_0)\leq -\lambda \int_m^{t_0} g(s)ds
\]
which clearly contradicts (\ref{eq:epsilon2}).

\textbf{Case 3 :} $C\neq 0$

Define $\hat{g}=g-C/\lambda$. One has
\[
\forall u,t \quad 0\leq u \leq t  \quad \hat{g}(t)-\hat{g}(u)=
g(t)-g(u)  \leq \int_u^t -\lambda g(s)+C ds = -\int_u^t \lambda \hat{g}(s) ds
\]
Thus $\hat{g}$ satisfies the conditions of Case 1 or Case 2, and as a
consequence
\[
\forall t\geq 0 \quad \hat{g}(t) \leq [\hat{g}(0)]^+e^{-\lambda t}
\]
The conclusion of the lemma follows by replacing $\hat{g}$ by
$g-C/\lambda$ in the above equation. $\Box$
